\documentstyle[12pt,twoside]{amsart}

\title
{Isomorphisms of jet schemes}
\author{Shihoko Ishii} 
\author{J\"org Winkelmann}

\address{Shihoko Ishii : Department of Mathematics, Tokyo Institute of
Technology, Oh-Okayama, Meguro, Tokyo, Japan
\newline
e-mail : shihoko@@math.titech.ac.jp}

\address{J\"org Winkelmann : Department of Mathematics, Universit\"at Bayreuth, Bayreuth, Germany
\newline
e-mail : joerg.winkelmann@@uni-bayreuth.de}

\thanks{The research of the first author was partially supported by Grant-in-Aid for Scientific Research (B) No.18340004 .  The research of the second author was  supported by  DFG Forschergruppe 790 ``Classification of algebraic
surfaces and compact complex manifolds''}

\newcommand{\bC}{{\Bbb C}}

\newcommand{\bN}{{\Bbb N}}
\newcommand{\bA}{{\Bbb A}}

\newcommand{\R}{{\cal R}}
\newcommand{\J}{{\cal J}}
\newcommand{\I}{{\cal I}}

\newcommand{\sTm}{{\spec K[t]/(t^{m+1})}}

\newcommand{\tm}{{k[t]/(t^{m+1})}}

\let \cedilla =\c
\renewcommand{\c}[0]{{\mathbb C}}  

\renewcommand{\o}[0]{{\mathcal O}} 

\newcommand{\spec}[0]{\operatorname{Spec}}

\def\to {\longrightarrow}

\newtheorem{thm}{Theorem}[section]

\newtheorem{lem}[thm]{Lemma}
\newtheorem{cor}[thm]{Corollary}
\newtheorem{prop}[thm]{Proposition}

\theoremstyle{definition}
\newtheorem{defn}[thm]{Definition}

\newtheorem{rem}[thm]{Remark}

\theoremstyle{remark}

\begin{document}

\maketitle

\markboth{\hfill SHIHOKO ISHII AND J\"ORG WINKELMANN \hfill}{\hfill
 ISOMORPHISMS OF JET SCHEMES  \hfill}

\begin{abstract}
If two schemes are isomorphic, then their $m$-jet schemes are isomorphic for all $m$. In this paper we consider the converse problem.
We prove that if an  isomorphism of the $m$-jet schemes is induced from a morphism of the base schemes, then the morphism of the base schemes is an isomorphism. 
But we also prove  that  just the existence of  isomorphisms between $m$-jet schemes does not yield the existence of an isomorphism between the base schemes.

\vskip 0.5\baselineskip

\noindent{R\'esum\'e} 
Si deux sch\'emas sont isomorphes, alors pour tout $m$ les sch\'emas de leurs $m$-jets
sont isomorphes. Dans cet article nous consid\'erons la question inverse.
Nous d\'emontrons que si un isomorphisme des sch\'emas de $m$-jets est induit
par un morphisme des sch\'emas de base, alors ce morphisme des sch\'emas de base est
un isomorphisme.
Mais nous d\'emontrons aussi que
l'existence d'un isomorphisme entre sch\'emas des $m$-jets 
n'implique pas l'existence d'un isomorphisme entre les sch\'emas 
de base. 

\vskip.3truecm
\noindent
AMS Mathematical Subject Classification 2000 : 14B05,  14J10, Keyword : jet schemes, arc spaces.

\end{abstract}

\section{Introduction}
\noindent
  From a given  scheme $X$ we obtain the $m$-jet scheme $X_m$ for every $m\in \bN\cup\{\infty\}$ (\cite{nash}). 
These jet schemes are considered to represent the nature of the scheme $X$ (see for example \cite{ein}, \cite{i6},  \cite{must01}).
It is clear that if we have an isomorphism $f: X\stackrel{\sim}\longrightarrow Y$ of schemes, then the induced morphism $f_m:X_m \to Y_m$ is also an isomorphism for every $m\in \bN\cup\{\infty\}$. 
Conversely, we consider the problem : if a morphism $f: X\longrightarrow Y$ induces an isomorphism 
$f_m:X_m \stackrel{\sim}\longrightarrow  Y_m$ for some $m\in \bN\cup\{\infty\}$, then is
$f$ an isomorphism ?
The answer is ``yes" and the proof follows immediately from the structure of jet schemes (see Corollary \ref{exist}).
Then, how about the problem without the assumption of the existence of a morphism $f$?
I.e., if there is an isomorphism $X_m\simeq Y_m$ for some $m$ or all 
$m\in \bN\cup\{\infty\}$, then is there an isomorphism $X\simeq Y$?
The answer is ``no".
We give an example of two non-isomorphic varieties $X$ and $Y$ with isomorphic $m$-jet schemes $X_m\simeq Y_m$ for every $m\in \bN\cup\{\infty\}$.

About jet schemes we use the terminology and notation in \cite{cr}.
In this paper a $k$-scheme is always a  scheme of finite type over an algebraically closed field $k$.
When we say variety, then it means an irreducible reduced separated scheme of finite type over $k$.

\section{Glossary of jet schemes}
\begin{defn}
  Let \( X \) be a $k$-scheme 
and $K\supset k$ a field extension.
  For \( m\in \bN \), a  \( k \)-morphism \( \sTm   \to X \) is called an \( m \)-{\it jet}
  of \( X \) 
 and a  \( k \)-morphism \( \spec K[[t]]\to X \) is called an  {\it arc} or $\infty$-{\it jet} of \( X \).
 For $m\in\bN\cup \{\infty\}$, the space of $m$-jets of  $X$ is denoted by $X_m$ (for details, see \cite{cr}).
\end{defn}

  Let \( G=\bA_k^1\setminus \{0\}=\spec k[s,s^{-1}] \) be the algebraic group with the multiplication. Then $X_m$ has a canonical action of $G$ for every 
  $m\in \bN\cup\{\infty\}$. 
  This action is induced from the ring homomorphism 
   \( \tm \to 
  k[s,s^{-1}, t]/(t^{m+1}) \), \( t\mapsto s t \).\vskip.5truecm
The following is well known by experts (see, for example \cite{ip}).

\begin{lem}
\label{cat}
For every $m\in \bN\cup\{\infty\}$, the  canonical projection 
$\pi_m:X_m\to X$ is the categorical quotient of the $m$-jet scheme $X_m$ by the action of $G$.

\end{lem}

For the definition of the categorical quotient, see \cite[Definition 0.5]{git}.

It is well known that for a $k$-scheme $X$, the $m$-jet scheme 
$X_m$ is affine over $X$ and has a graded structure, i.e, if we write 
$X_m=\spec {\R^{(m)}}$, then $\R^{(m)}$ has a structure of graded $\o_X$-algebra
$\R^{(m)}=\oplus_{i\geq 0} \R_i^{(m)}$  with $\R_0=\o_X$.
The truncation morphism $\psi_{m, m'}:X_{m} \to X_{m'} $  $(m'<m)$ corresponds to a homomorphism of graded algebras 
$\psi^*_{m,m'}:\R^{(m')}\to \R^{(m)}$.
In particular, by $\psi^*_{m,m-1}$, $\R^{(m)}$ is an $\R^{(m-1)}$-algebra and generated by $m$-th homogeneous part $\R_m^{(m)}$ over $\R^{(m-1)}$.
Therefore, the restriction $\R_i^{(m-1)}\to \R_i^{(m)}$ of $\psi^*_{m,m-1}$ 
is surjective for $i\leq m-1$.

\begin{lem}
\label{exact}
  Let $X$ be a  $k$-scheme. 
  Under the  notation above, let $\J^{(m)}\subset \R^{(m)}$ be the $\R^{(m)}$-ideal 
  defining $\psi_{m,m-1}^{-1}(\sigma_{m-1}(X))$, where $\psi_{m,m-1}:
 X_m\to X_{m-1}$ is the truncation morphism and $\sigma_{m-1}:X\to X_{m-1}$ 
 is the canonical section of the projection $\pi_{m-1}:X_{m-1}\to X$. 
  Let $\J^{(m)}=\oplus _{i\geq 0}\J_i^{(m)}$ be the homogeneous decomposition.
    Then, $\J^{(m)}$ is generated by  $\oplus _{i=1}^{m-1}
 \R^{(m-1)}_i$ and 
there exists an exact sequence of $\o_X$-modules:
  $$0\to \J_m^{(m)} \to \R_m^{(m)} \to \Omega_{X/k} \to 0.$$
  \end{lem}
  
  \begin{pf}

  The closed subscheme $\sigma_{m-1}(X)\subset X_{m-1}$ is defined  
  by the ideal $\I=\oplus_{i\geq 1} \R^{(m-1)}_i $.
  Here, we note that $\I$ is generated by $\oplus _{i=1}^{m-1}
 \R^{(m-1)}_i$, since $\R^{(m-1)}$ is generated by $\oplus _{i=1}^{m-1}
 \R^{(m-1)}_i$ as an $\o_X$-algebra.
  Therefore, we obtain that 
 $\J^{(m)}$ is generated by  $\oplus _{i=1}^{m-1}
 \R^{(m-1)}_i$.

For the next statement, note that there is a canonical isomorphism
 $$\psi^{-1}_{m,m-1}(\sigma(X)) \simeq \spec S(\Omega_{X/k}),$$ 
 where $S(\Omega_{X/k})$ is the 
 symmetric $\o_X$-algebra defined by $\Omega_{X/k}$
(see for example \cite{ip}, where the proof is similar to that in \cite[Example 2.5]{em}).
Then, it follows an isomorphism  
 $\R^{(m)}/\J^{(m)}\simeq S(\Omega_{X/k})$ of graded 
 $\o_X$-algebras.
 Here we note that the isomorphism sends the part of  degree $mi$ of 
 $\R^{(m)}/\J^{(m)}$ to the part of degree $i$ of $S(\Omega_{X/k})$.
 By taking the generating parts of the both graded algebras, we have 
 $$\R_m^{(m)}/\J_m^{(m)}\simeq \Omega_{X/k}.$$
 \end{pf}

\begin{prop}
\label{mainprop}
Let $X$ be an affine variety with $\Omega_{X/k}\simeq \o_X^{\oplus n}.$
Then, $X_m\simeq X\times \bA_k^{mn}$ for every $m\in \bN$.
\end{prop}

\begin{pf}
As $X_1\simeq \spec S(\Omega_{X/k})$, it follows that $X_1\simeq X\times \bA_k^n$.
For $m\geq 2$, denote $X_m$ by $\spec \R^{(m)}$ and the $\R^{(m)}$-ideal defined in 
Lemma \ref{exact} by $\J^{(m)}$.
Then, by Lemma \ref{exact}, we have the exact sequence:
\begin{equation}\label{ex}
0\to \J^{(2)}_2\stackrel{\iota}\longrightarrow \R^{(2)}_2\stackrel{\psi}\longrightarrow \Omega_{X/k}\to 0.
\end{equation}
Since $\Omega_{X/k}\simeq \o_X^{\oplus n}$ and $X$ is affine,
there is a section $ s:\Omega_{X/k}\to \R^{(2)}_2$ of $\psi$ and therefore the exact sequence (\ref{ex}) splits.
Since $X$ is a non-singular variety, the homomorphism $\R_i^{(1)}\to \R_i^{(2)}$ $(i=0,1)$ is  injective and therefore bijective by the note before Lemma \ref{exact}, it follows that $\R_i^{(2)}=\R_i^{(1)}$ for $i=0,1$ with identifying by the bijection.
Then, by Lemma \ref{exact} $$\J_2^{(2)}=\R_1^{(1)}\cdot \R_1^{(2)}+\R_2^{(1)}\cdot \R_0^{(2)}$$
$$= 
\R_1^{(1)}\cdot \R_1^{(1)}+\R_2^{(1)}\cdot \R_0^{(1)}\subset \R^{(1)}.$$ 
Therefore, by the splitting exact sequence (\ref{ex}), an element of  $\R_2^{(2)}$ is the sum of an element of  $\R^{(1)}$ and an element of $s(\Omega_{X/k})$ which is globally generated by 
$n$ elements over $\o_X$.
As $\R^{(2)}$ is generated by $\R_2^{(2)}$ over $\R^{(1)}$,
we have a surjection 
$$\R^{(1)}[\theta_1, \ldots, \theta_n] \to \R^{(2)},$$
where $\theta_1, \ldots ,\theta_n$ are indeterminates.
Thus, we obtain a closed immersion $X_2\hookrightarrow X_1\times \bA^n$.
Considering the dimension of the both varieties, we have 
$X_2\simeq  X_1\times \bA^n$.
For $m>2$, similarly we obtain $X_m\simeq X_{m-1}\times \bA^n$, which shows $X_m\simeq X\times \bA_\bC^{nm}$ for every $m$.
\end{pf}

\begin{rem}
This proposition can be proved by using the fact that $X_m$ has a $X_1\times
_X{X_{m-1}}$-torsor structure over $X_{m-1}$.
But we prefer  our present proof, because it is elementary.
\end{rem}
 
 \section{Isomorphism problems}
 \noindent
  A $G$-equivariant isomorphism yields an isomorphism of 
   the categorical quotients. 
  Therefore, by Lemma \ref{cat} we obtain:

 \begin{prop}
 \label{equiv}
 Let $X$ and $Y$ be two schemes over $k$. 
 If there exists a $G$-equivariant isomorphism $X_m\stackrel{\sim}\longrightarrow Y_m$ of $m$-jet schemes for some $m\in \bN\cup\{\infty\}$,
 then there is an isomorphism  $X\stackrel{\sim}\longrightarrow Y$.
  \end{prop}
 
 The induced morphism  $f_m:X_m \to Y_m$ from a morphism $f:X\to Y$ is $G$-equivariant.
 Therefore, by the previous proposition and the universality of the categorical quotient, 
 we obtain the following:

 \begin{cor}
 \label{exist}
 Let $f:X\to Y$ be a morphism of schemes over $k$. If the induced morphism $f_m:X_m\to Y_m$ is isomorphic for some $m\in \bN\cup\{\infty\}$,
 then the morphism $f$ is an isomorphism.
  \end{cor}

\begin{rem}
This corollary can be proved directly by using the fact that the morphism of the base spaces induces the morphism of the sections in the jet-schemes.
\end{rem} 

Now let us be just given an isomorphism of $m$-jet schemes and we consider if it induces an isomorphism of base schemes.
The following is a counter example for this problem.
We use the counter example of the cancellation  problem called Danielewski's example.

\begin{thm} Let $X$ and $Y$ be hypersurfaces in $\bA_{\bC}^3$ defined by 
$ xz-y^2+1=0 $ and $x^2z-y^2+1=0$, respectively.
Then, $X\not\simeq Y$ but $X_m\simeq Y_m$ for every $m\in \bN\cup\{\infty\}$.
\end{thm}

\begin{pf}
  By the work of Danielewski it is known that $X,Y$ are non-singular, $X\not\simeq Y$ and $X\times \bA_{\bC}^1\simeq Y\times \bA_{\bC}^1$ (see for example, \cite{ml}, \cite{sy}).
  Therefore, we have only to prove that
  $X_m\simeq X\times \bA_{\bC}^{2m}, \ \ \ Y_m\simeq Y\times \bA_{\bC}^{2m}$
  for $m<\infty$ and
  $X_\infty\simeq Y_\infty.$ 
   According to  Danielewski's idea, define  actions of the additive group $\bA_{\bC}^1$ on $X$ and $Y$ as follows:
$$\bA_{\bC}^1\times X\to X,\ \ \ (t, x,y,z)\mapsto (x, y+xt,z+2yt+xt^2),$$
$$\bA_{\bC}^1\times Y\to Y,\ \ \ (t, x,y,z)\mapsto (x, y+x^2t,z+2yt+x^2t^2).$$
Then the  actions are free and the number of the orbits for fixed $x\neq 0$ is one and for $x=0$ it is two for both actions.
By this we have principal fiber bundles 
$\varphi_1:X\to Z,\ \ \ \varphi_2:Y\to Z,$
where $Z=\bA_{\bC}^1\cup\bA_{\bC}^1$ is the  line with bug-eyes; i.e., the union of 
the two lines $\bA_{\bC}^1$ with  patching by $id:\bA_{\bC}^1\setminus \{0\}\simeq \bA_{\bC}^1\setminus \{0\}$.
As the morphisms $\varphi_i$ $(i=1,2)$ are smooth, 
we have surjections of tangent sheaves:
$$T_{X/\bC}\to \varphi_1^*T_{Z/\bC}, \ \ \ T_{Y/\bC}\to \varphi_2^*T_{Z/\bC}.  $$
Let $L_i$ $(i=1,2)$ be their kernels, respectively. 
Since $T_{Z/\bC}\simeq \o_Z$ and $X$, $Y$  are  affine varieties, we have the splittings:
$$T_{X/\bC}\simeq \o_X\oplus L_1,\ \ \ T_{Y/\bC}\simeq \o_Y\oplus L_2.$$
By this, we obtain the canonical sheaves $\omega_X\simeq L_1^{-1}$,
$\omega_Y\simeq L_2^{-1}$.
On the other hand, it is known that 
$\omega_{\bA_{\bC}^3}\simeq\o_{\bA_{\bC}^3}$, $\o_{\bA_{\bC}^3}(X)\simeq \o_{\bA_{\bC}^3}$ and $\o_{\bA_{\bC}^3}(Y)\simeq \o_{\bA_{\bC}^3}$.
Then, 
by Adjunction Formula, we obtain
 $$\omega_X\simeq\omega_{\bA_{\bC}^3}(X)\otimes\o_X\simeq \o_X, \ \ 
\omega_Y\simeq\omega_{\bA_{\bC}^3}(X)\otimes\o_Y\simeq \o_Y,$$
which shows that $L_i$'s are trivial.
Hence, we have $\Omega_{X/\bC}\simeq \o_X^{\oplus 2}$ and 
$\Omega_{Y/\bC}\simeq \o_Y^{\oplus 2}$.
By Proposition \ref{mainprop}, we have 
$X_m\simeq X\times \bA_\bC^{2m}, Y_m\simeq Y\times \bA_\bC^{2m}$ for every $m\in \bN$. 
For the proof of $X_\infty\simeq Y_\infty$, we have only to take the projective limit of the isomorphisms $X_m\simeq Y_m$ which are compatible with the truncation morphisms.

\end{pf}

\section*{Acknowledgements}
The basic idea of our result came up when the authors and Professor Laszlo Lempert were talking at Hayama Conference.
We would like to thank Professor Lempert  and the organizers of the conference 
for giving us a chance to collaborate on this topic.
The first author is grateful to Professor Teruo Asanuma for providing with references on the cancellation problem.

\end{document}